\documentclass[11pt]{article}
\usepackage{amsmath}
\usepackage{amssymb}
\usepackage{amsfonts}
\usepackage{eucal}
\usepackage[usenames]{color}
\usepackage{graphicx}
\numberwithin{equation}{section}
\newtheorem{Theorem}{Theorem}[section]
\newtheorem{Proposition}[Theorem]{Proposition}

\newtheorem{Lemma}[Theorem]{Lemma}
\newtheorem{Remark}[Theorem]{Remark}

\begin{document}
\title{On the representation of an integrated  Gauss-Markov process}
\author{Mario Abundo\thanks{Dipartimento di Matematica, Universit\`a  ``Tor Vergata'', via della Ricerca Scientifica, I-00133 Rome, Italy.
E-mail: \tt{abundo@mat.uniroma2.it}}
}
\date{}
\maketitle

\begin{abstract}
\noindent We find a representation of the integral of a
Gauss-Markov process in the interval $[0,t],$ in terms of Brownian
motion. Moreover, some connections with first-passage-time
problems are discussed, and some examples are reported.

\end{abstract}


\section{Introduction}
In this short note, we consider a real continuous Gauss-Markov
process $X(t)$ of the form:
\begin{equation}
X(t) = m(t) + h_2(t) B(\rho (t)), \ t \ge 0
\end{equation}
where: \par\noindent
$\bullet \ B(t)$ is a standard Brownian motion (BM); \par\noindent
$\bullet \ m(t)= E(X(t))$ is continuous for every $t\ge 0;$ \par\noindent
$ \bullet \ $ the covariance $c(s,t) := E [(X(s) -m(s)) (X(t)-m(t))]$
is continuous for every $0 \le s <t ,$ with \par\noindent \ \ \   $c(s,t)=h_1(s) h_2(t); $ \par\noindent
$ \bullet \  \rho(t)= h_1(t)/h_2(t) $ is a monotonically increasing function and $h_1(t) h_2(t) >0 .$ \par\noindent
Notice that a special case of Gauss-Markov process is the Ornstein-Uhlenbeck (OU) process, and in fact
any Gauss-Markov process can be represented in terms of a OU process (see e.g. \cite{ric:smj08}).  \par\noindent
Our aim is to find a representation of
\begin{equation}
Y(T): = \int _ 0 ^T X(s) ds, \ T >0 ,
\end{equation}
in terms of Brownian motion.
Notice that the integrated process $Y(T)$ is equal to $\overline X _T \cdot T,$ where $\overline X _T$ is the time average of $X(t)$ in the interval $[0,T].$
\par\noindent
The study of $Y(T)$  has interesting applications in Biology, for
instance in the framework of diffusion models for neural activity;
if one identifies $X(t)$ with the neuron voltage at time $t,$
then, $Y(T) / T $  represents the time average of the neural
voltage in the interval $[0,T].$ Another application can be found
in Queueing Theory, if  $X(t)$ represents the length of a queue at
time $t;$ then, $Y(T)$ represents the cumulative waiting time
experienced by all the ``users'' till the time $T.$
\par\noindent
As for an example from Economics, let us suppose that the
variable $t$ represents the quantity of a commodity that producers
have available for sale, then $Y(T)$
provides a measure of the total value that
consumers receive from consuming the amount $T$  of the product. \par
Really, in certain applications, it is interesting to study first-passage time (FPT) problems for
the integrated process $Y(T);$ to this end,
it is useful to dispose of an explicit representation of $Y(T).$ The results of this paper generalize those of \cite{abundo:ija08}.

\section{Main Results}
We begin with stating and proving the following:

\begin{Lemma} \label{lemma}
Let $f(t)$ a continuous bounded deterministic function, with $f(t) \neq 0$ for every $t \ge 0,$  then
\begin{equation}
I(t) := \int _0 ^t f(s) B(s) ds
\end{equation}
is normally distributed with  mean zero and variance $\gamma (t),$ where
$\gamma (t)= \int _0 ^t (R(t) - R(s)) ^2 ds $ and $R(t)= \int _0 ^t f(s) ds .$ Moreover, if $\gamma (+ \infty ) = + \infty, $
then there exists a BM $\widetilde B (t)$ such that  $I(t) = \widetilde B ( \gamma (t) ).$
\end{Lemma}

\noindent{\it Proof.} \  We observe that $I(t)$ is a Gaussian process with zero mean and variance
$$ V(t) :=  Var (I(t)) = Cov \left ( \int _0 ^t f(s) B(s) ds, \int _0 ^t f(u) B(u) du \right ) $$
$$ = E \left ( \int _0 ^t f(s) B(s) ds \ \cdot \  \int _0 ^t f(u) B(u) du \right )
  = \int _0 ^t ds \int _0 ^t du E( f(s)B(s) f(u) B(u)) . $$
Since $E( f(s)B(s) f(u) B(u)) = f(s)f(u) \min (s,u),$ we get:
$$ V(t)= \int \int _ { \Delta _1 } f(s) f(u) \ u \ ds du + \int \int _ { \Delta _2 } f(s) f(u) \ s  \ ds du  ,$$
where $\Delta _1 = \{ (s,u) \in [0, + \infty ) \times [0, + \infty )  : 0 \le s \le t, \ 0 \le u \le s \}$ and
$\Delta _2 = \{ (s,u) \in [0, + \infty ) \times [0, + \infty ) : 0 \le s \le t,  u \ge s \}.$ Thus, by calculation, we
obtain:
$$ V(t)= 2 \int _ 0 ^t f(s) ds \int _0 ^s f(u) \ u \  du .$$
As easily seen,  $V(t)$ and $\gamma (t)$ have the same derivative,
so the equality  $V(t)= \gamma (t)$ follows for any $t \ge 0,$  since $V(0)=\gamma (0)=0.$ \par\noindent
Let $T>0$ fixed;
by using ${\rm It \hat o}$'s formula we get:
$$ I(T) = \int _0 ^T f(s) B(s) ds = R(T) B(T) - \int _0 ^T R(s) dB(s) = \int _0 ^T (R(T) - R(s) ) dB(s) .$$
For $t \le T,$ let us consider now the continuous martingale $M_t$ having differential $dM_t = (R(T) -R(t)) dB_t;$
for $s \le t \le T$ we have:
\begin{equation} \label{ineq}
\gamma (t)= \int _0 ^t (R(t) -R(s))^2 ds \le \int _0 ^T (R(T) -R(s))^2 ds = \langle M \rangle_T,
\end{equation}
where $\langle M \rangle_t$ denotes the quadratic variation of
$M_t.$ Indeed, the inequality \eqref{ineq} easily follows from the
fact that, since $R'(t)=f(t) \neq 0,$ the function $R(t)$ is
monotone (increasing or decreasing), so $T \ge t$ implies $(R(T)
-R(s))^2 \ge (R(t) - R(s))^2.$ \par\noindent If $\gamma ( + \infty
) = + \infty,$ from inequality \eqref{ineq} we obtain $\langle M
\rangle_ \infty = + \infty ;$ then, by the Dambis, Dubins-Schwarz
Theorem (see e.g. \cite{revuzyor:con91}) there exists a BM
$\widetilde B $ such that $M_t = \widetilde B ( \langle M \rangle_
t).$ Thus, $I(T) = M_T = \widetilde B ( \langle M \rangle_ T);$
finally, taking $t=T,$ we obtain $I(t) = M_t = \widetilde B (
\langle M \rangle_ t) =  \widetilde B ( \gamma (t)).$

\hfill $\Box$
\bigskip

\noindent As a corollary of the previous lemma, we obtain our main result:
\begin{Proposition} \label{mainproposition}
Let $X(t)$ be a Gauss-Markov process given by (1.1), and suppose that $h_1, \ h_2$ are continuous,  $\rho$  is a differentiable increasing function; then $Y(t)= \int _ 0 ^t X(s) ds $
is normally distributed with  mean $M(t)= \int _0 ^t m(s) ds$  and variance $\gamma _1 (\rho(t)),$ where
$\gamma _1 (t)= \int _0 ^t (R_1(t) - R_1(s)) ^2 ds $ and $R_1 (t)= \int _0 ^t h_2(\rho ^{-1} (s))/ \rho '( \rho ^{-1} (s)) ds .$
Moreover, if $\gamma _1 (+ \infty ) = + \infty, $
then there exists a BM $\widehat B (t)$ such that  $Y(t) = M(t) + \widehat B ( \gamma _1 (\rho(t)) ).$ Thus, $Y(t)$ is still
Gauss-Markov.
\end{Proposition}

\noindent{\it Proof.} \ We have:
$$ Y(t)= \int _0 ^t X(s) ds $$
$$ = \int _0 ^t m(s) ds + \int _0 ^t h_2
(s) B ( \rho (s)) ds = M(t) + \int _0 ^ { \rho (t)} h_2 ( \rho
^{-1} (s))/ \rho '( \rho ^{-1} (s)) B(s) ds ,$$
where we have used a variable change in the integral.
Then, the proof follows by using  Lemma \ref{lemma} with $f(t)= h_2 ( \rho
^{-1} (t))/ \rho '( \rho ^{-1} (t)).$

\hfill $\Box$
\bigskip

\begin{Remark}  If we consider
the time average $\overline X_T = \frac 1 T \int _0 ^T X(s) ds ,$ by  Proposition \ref{mainproposition} we get
$\overline X_T = \frac 1 T Y(T) = \frac 1 T \left [ M(T) + \widehat B (\gamma _1 (T)) \right ],$ namely
$\overline X_T $ is normally distributed with mean $M(T)/T$ and variance $\gamma _1 (T) / T^2 .$  In particular,
if $X(t)$ is BM, one obtains $\overline X_T \sim {\cal N } (0,T/ 3)$ (cf.  \cite{abundo:ija08}).
\end{Remark}

\begin{Remark} Notice that, if $\gamma _1 ( \infty ) = \infty, $ than the FPT of $Y(t)$
over a continuous boundary $S(t)> 0,$ i.e.
$ \tau _S = \inf \{ t >0 : Y(t) \ge S(t) \},$ is nothing but the FPT of
$\widehat B ( \gamma _1 (t) )$ over $ \bar S(t)= S(t) - M(t) ,$ or equivalently  $\gamma _1 (\tau _S)=
\inf \{ u >0 : \widehat B  (u) > \bar S ( \gamma _1 ^{-1} (u))\} .$
\end{Remark}

\section{A Few Examples}

\noindent {\bf Example 1} (Brownian motion with drift)
\par\noindent Let be $X(t)= \mu t  + B(t),$  then $m(t)= \mu t , \
h_1(t) = t,  \ h_2 (t)=1$ and $\rho(t)=t.$ Moreover, $R_1(t)=\int
_0^t ds =t$ and $\gamma _1 (t) = \int _0 ^t (t-s)^2 ds = t^3/3 .$
Thus, $Y(t) = \mu t^2 /2 + \widehat B (t^3/ 3)$ (cf.
\cite{abundo:ija08}).
\bigskip

\noindent {\bf Example 2} (Ornstein-Uhlenbeck process)\par\noindent
Let $X(t)$ be the solution of the SDE:
$$ dX(t)= - \mu (X(t)-\beta)dt + \sigma dB_t, \ X(0)=x $$
where $\mu , \sigma >0$ and $\beta \in (- \infty, + \infty ).$
The explicit solution is (see e.g.  \cite{abundo:stapro12} ):
$$X(t)=  \beta + e^{ - \mu t } [ x-\beta +  \widetilde B(\rho (t)] $$
where $\widetilde B$ is Brownian motion and $\rho (t)= \frac {\sigma ^2} {2 \mu} \left (e ^{2 \mu
t } -1 \right ) .$ So, $X(t)$ is a Gauss-Markov process with
$m(t)=  \beta + e^{ - \mu t } (x- \beta ), \ h_1(t)= \frac {\sigma ^2} {2 \mu} \left (e ^{\mu
t } -e^{- \mu t } \right ), \
h_2(t)= e^{ - \mu t }$ and $c(s,t)= h_1(s) h_2(t).$ By calculation, we obtain:
$$M(t)= \int _0 ^t \left ( \beta + e^{ - \mu s } (x- \beta ) \right ) \ ds= \beta t + \frac {(x - \beta)} \mu  \left (1 - e^{ - \mu t } \right ) , $$
$$R_1(t) = \int _0 ^t  e^{- \mu \rho ^{-1} (s)} ( \rho ^{-1})'(s) ds = \frac {1- e ^{- \nu \rho ^{-1} (t) }  } { \mu} , $$
$$ \rho ^{-1} (s)= \frac 1 { 2 \mu} \ln \left (1 + \frac{ 2 \mu} {\sigma ^2} s \right ) , $$
$$ \gamma_1(t)= \frac 1 { \mu ^2} \int _0 ^t \left ( e^{- \mu \rho ^{-1} (t) } - e^{ - \mu \rho ^{-1} (s)} \right ) ^2 ds =
\frac 1 { \mu ^2} \int _0 ^t \left ( \frac 1 { \sqrt { 1+ 2 \mu t / \sigma ^2 } } - \frac 1 { \sqrt { 1+ 2 \mu s / \sigma ^2 } } \right ) ^2 ds $$
$$=  \frac {\sigma ^2 t } {\mu ^2 (\sigma ^2 + 2 \mu t ) } - \frac {2 \sigma ^2 } {\mu ^3 \sqrt { 1+ 2 \mu t / \sigma ^2 } }
\left ( \sqrt { 1+ 2 \mu t / \sigma ^2 } -1 \right ) + \frac { \sigma ^2 } {2 \mu ^3 } \ln \left ( 1+ 2 \mu t / \sigma ^2 \right ) .$$
Then, by Proposition \ref{mainproposition}, we get that
$Y(t)= \int _0 ^t X(s) ds $ is normally distributed with mean
$M(t)$
and variance $\gamma _1 (\rho (t))  .$ Moreover, since
$\lim _ { t \rightarrow + \infty} \gamma _1(t) = + \infty ,$
there exists a BM $\widehat B(t) $ such that $Y(t)= M(t) + \widehat B \left ( \gamma _1 (\rho (t)) \right ).$
\bigskip

\noindent{\bf Example 3} (Brownian bridge) \par\noindent
For $T>0$ and given $a, b,$ let $X(t)$ be the solution of the SDE:
$$ dX(t)= \frac {b - X(t)} {T-t } \ dt +  dB_t, \ 0 \le t \le T, \ X(0)=a .$$
This is a transformed BM with fixed values at each end of the interval $[0,T], \ X(0)=a$ and $X(T)=b.$
The explicit solution is (see e.g. \cite{revuzyor:con91}):
$$ X(t) =a \left ( 1 - t/T \right ) + bt /T + (T-t) \int _0 ^t \frac 1 { T-s} dB_s $$
$$= a \left ( 1 - t/T \right ) + bt /T + (T-t) \widetilde B \left ( \frac { t} {T(T-t) } \right )  , \ 0 \le t \le T ,$$
where $\widetilde B$ is BM.
So, $X(t)$ is a Gauss-Markov process with:
$$m(t)= a \left ( 1 - t/T \right ) + bt /T, \ c(s,t)= h_1(s)h_2(t) \ {\rm  with} \  h_1(t) = t/T , \ h_2(t) =T-t, \ \rho (t) = \frac t { T(T-t)} .$$
By calculation, we obtain:
$$M(t)= at + \frac {b-a } {2T }  t^2, $$
$$R_1(t) = \frac {T^3t (2+Tt) } {2(1+Tt)^2}  , $$
$$ \rho ^{-1} (s)= \frac {T^2 s } {1+Ts } , $$
$$ \gamma_1(t)= \int _0^t \left ( \frac {T^3t (2+Tt) } {2(1+Tt)^2}- \frac {T^3s (2+Ts) } {2(1+Ts)^2} \right )^2   ds.$$
Then, by Proposition \ref{mainproposition}, we get that $Y(t)=
\int _0 ^t X(s) ds $ is normally distributed with mean $M(t)$ and
variance $\gamma _1 (\rho (t))  .$ By a straightforward, but
boring calculation, it can be verified that  $\lim _ { t
\rightarrow + \infty} \gamma _1(t) = + \infty ,$ so there exists a
BM $\widehat B(t) $ such that $Y(t)= M(t) + \widehat B \left (
\gamma _1 (\rho (t)) \right ).$
\bigskip

\noindent{\bf Example 4} (Generalized Gauss-Markov process) \par\noindent
Let us consider the diffusion $X(t)$ which is the solution of the SDE:
$$ dX(t)= m'(t) dt + \sigma (X(t)) dW_t , \ X(0) = m(0) $$
where $W_t$ is BM and $\sigma (x)>0$ is a smooth deterministic function. We suppose that
$\rho (t)= \langle X \rangle _t = \int _0 ^t \sigma ^2 (X(s)) ds , $ i.e. the quadratic variation of $X(t),$ is
 increasing to $\rho ( + \infty)= + \infty .$
By using the Dambis, Dubins-Schwarz Theorem, it follows that
$X(t)= m(t)+ B( \rho (t)), \ t \ge 0 ,$ where $\rho(t)$ is not necessarily deterministic, but it can be a random function.
For this reason, we call $X(t)$ a generalized Gauss-Markov process.
Denote by $A$ the ``inverse'' of the random function $\rho ,$ that is,
$ A(t) = \inf \{ s >0 : \rho (s) >t \};$ since $\rho(t)$ admits derivative and $\rho'(t)=\sigma ^2 (X(t)) > 0,$ also $A'(t)$ exists and
$ A'(t)= 1/ \sigma ^2 (X(A(t)));$
we focus on the case when
there exist  deterministic continuous  functions $\alpha (t), \ \beta (t) $ (with  $\alpha(0)= \beta (0))$ and
$\alpha _1(t), \ \beta _1 (t),$ such that, for every $t \ge 0:$
$$ \alpha (t), \  \beta (t) \ {\rm are \ increasing ,} \ \alpha (t) \le \rho (t) \le \beta (t) , \ {\rm and} \
\alpha _1 (t) < A'(t) < \beta _1  (t).$$
Since $\rho (t)$ is not deterministic, we cannot obtain exactly the distribution of  $ \int _0 ^t X(s) ds,$ however we are able
to find bounds to it. In fact, we have:
$$ \int _0^t X(s) ds = \int _0 ^t m(s) ds + \int _0 ^t B(\rho (s)) ds= \int _0 ^t m(s) ds + \int _0 ^ {\rho (t)} B(v) A' (v) dv $$
We can use the arguments of Lemma \ref{lemma}  with $f(v)=A' (v), \
R_1(t)= \int _0 ^t A' (s) ds ,$ and $\gamma _1 (t)= \int _0 ^t (R_1(t)- R_1(s))^2 ds;$  by assumptions we get
$\int _ 0 ^t \alpha _1 (s) ds \le R_1(t) \le \int _ 0 ^t \beta _1 (s) ds  .$
Thus, we conclude that
$ \int _0 ^t X(s) ds$ is normally distributed with mean $ M(t)= \int _0 ^t m(s) ds$ and variance $\gamma _1 (\rho(t)) ,$
which is bounded between $\gamma _1(\alpha (t))$ and $\gamma _1(\beta (t)).$ The closer $\alpha (t)$ to $\beta (t),$ the better the approximation above;
for instance, if $\sigma (x)= 1+ \epsilon \cos ^2 (x), \  \epsilon >0 ,$ we have $\rho(t) = \int _0 ^t (1+  \epsilon \cos ^2 (X(s)))^2 ds$ and so
$\alpha (t)= t, \ \beta (t)= (1+ \epsilon)^2 t , \ \alpha _1(t) = 1 / (1+ \epsilon )^2, \ \beta _1(t)= 1. $
The smaller is $\epsilon,$ the closer $\gamma_1 (\alpha (t))$ to $\gamma _1 (\beta (t)).$


\begin{thebibliography}{spc}


\bibitem [1] {abundo:stapro12}
Abundo, M., 2012. \newblock
An inverse first-passage problem for one-dimensional diffusions with random
      starting point.
\newblock{Statistics  and Probability Letters} 82 (1), 7--14.

\bibitem [2] {abundo:ija08}
Abundo, M., 2008. \newblock
On the distribution of the time average of a  jump-diffusion process.
\newblock{International Journal of Applied  Mathematics (IJAM)} 21 (3), 447--454.


\bibitem[3]
{revuzyor:con91}
Revuz, D. and Yor, M., 1991. \newblock
Continous martingales and Brownian motion. \newblock
Springer-Verlag, Berlin Heidelberg.

\bibitem[4]
{ric:smj08}
Nobile, A.G., Pirozzi, E., Ricciardi, L.M. , 2008. \newblock
Asymptotics and evaluations of FPT densities through varying boundaries for Gauss-Markov processes. \newblock
{Scientiae Mathematicae Japonicae} 67, (2), 241--266.






\end{thebibliography}
\end{document}